\documentclass[preprint,10pt,authoryear]{elsarticle}
\usepackage{graphicx, verbatim, url, amssymb, amsmath, amsfonts}
\usepackage{algorithm}
\usepackage{geometry}
\newcommand{\fv}{\mbox{\boldmath  $v$}}
\newcommand{\fu}{\mbox{\boldmath  $u$}}
\newcommand{\fw}{\mbox{\boldmath  $w$}}
\newcommand{\fy}{\mbox{\boldmath  $y$}}
\newcommand{\fr}{\mbox{\boldmath  $r$}}
\newcommand{\fx}{\mbox{\boldmath  $x$}}
\newcommand{\ww}{\sc internet \rm}

\newcommand{\ds}{\displaystyle}
\renewcommand{\t}{^T}
\newcommand{\prot}{{\tt protein }}
\newcommand{\inte}{{\tt internet\ }}
\newcommand{\www}{{\tt www\ }}
\newcommand{\dblp}{{\tt dblp\ }}

\journal{Internat Mathematics}

\begin{document}
\begin{frontmatter}

\title{Computing the smallest eigenpairs of the graph Laplacian}

\author{Luca Bergamaschi}
\address{Department of Civil Environmental and  Architectural Engineering , University of Padua, Italy\\
{\tt e-mail} luca.bergamaschi@unipd.it}

\author{Enrico Bozzo}
\address{Department of Mathematics and Computer Science, University of Udine, Italy\\
{\tt e-mail} enrico.bozzo@uniud.it}

\author{Massimo Franceschet}
\address{Department of Mathematics and Computer Science, University of Udine, Italy \\
{\tt e-mail} massimo.franceschet@uniud.it}

\begin{abstract}
The graph Laplacian, a typical representation of a network, is an important matrix that can tell us much about the network structure. 
In particular its eigenpairs (eigenvalues and eigenvectors) incubate precious topological information about the network at hand, 
including connectivity, partitioning, node distance and centrality. Real networks might be very large in number of nodes (actors); 
luckily, most real networks are sparse, meaning that the number of edges (binary connections among actors) are few with respect to 
the maximum number of possible edges. In this paper we experimentally compare three state-of-the-art algorithms for computation of 
a few among the smallest eigenpairs of large and sparse matrices: the Implicitly Restarted Lanczos Method, which is the current implementation 
in the most popular scientific computing environments (Matlab $\slash$ R), the Jacobi-Davidson method, and the 
Deflation Accelerated Conjugate Gradient method. We implemented the algorithms in a uniform programming setting and tested them over 
diverse real-world networks including  biological, technological, information, and social networks. It turns out that the Jacobi-Davidson 
method displays the best performance in terms of number of matrix-vector products and CPU time. 
\end{abstract}


\end{frontmatter}

\section{Introduction} \label{introduction}

The bi-directional link between the relatively new discipline of network science and the well-consolidated field of matrix algebra 
is intriguing and promising \citep{M11,BH12}. The big challenge is to bridge network science and matrix algebra in a synergy. 
Can we apply results and methods of matrix algebra to investigate the properties of networks? Do real networks, with their 
universal architectures, represent a class of algebraic structures (matrices) for which results and methods of matrix algebra 
can be improved or specialized? Clearly, both network science and matrix algebra would benefit from this synergistic approach. 
Network science would gain additional insight in the structure of real networks, while matrix algebra would obtain more challenging applications.

Networks, in their basic form of graphs of nodes and edges, can be represented as matrices. 
The most common representation of a graph consists of the graph adjacency matrix, where the entries of the 
matrix that are not null represent the edges of the graph. Often, it is convenient to represent a graph with its Laplacian matrix, 
which places on the diagonal the degrees of the graph nodes (the number of connections of the nodes) and elsewhere information about the 
distribution of edges among nodes in the graph. The Laplacian matrix, and in particular The Laplacian matrix, and in particular its 
smallest eigenpairs (eigenpairs relative to the smallest eigenvalues), turn up in many different places in network science. 
Examples include random walks on networks, resistor networks, resistance distance on networks, current-flow closeness and betweenness centrality 
measures, graph partitioning, and network connectivity \citep{GBS08,BF2012,N10}. 

Real networks might be very large; however, they are typically also very sparse. Moreover, generally, 
not the entire matrix spectrum is necessary, but only a few eigenpairs, either the lowest of the largest, are enough. 
A number of iterative procedures, based on a generalization of the well-known power method, have been recently developed  to 
compute a few eigenpairs of a large and sparse matrix. 

In this paper, we experimentally analyze three important iterative methods: 
(i) the Implicitly Restarted Lanczos Method, 
(ii) the Jacobi-Davidson method, and 
(iii) the Deflation Accelerated Conjugate Gradient method. 
We implement these methods in a uniform programming environment and experimentally compare them on four Laplacian matrices of networks 
arising from realistic applications. The real networks include a biological network (a protein-protein interaction network of yeast), 
a technological network (a snapshot of the Internet), an information network (a fragment of the Web), 
and a social network (the whole collaboration network among Computer Science scholars).

%We selected to use for the three methods the established implementations without making  optimization
%to  any of them. However, it is worth mentioning that much work is being devoted particularly  to the Arnoldi method (the non symmetric
%counterpart of the Lanczos Method) in order to reduce its computational cost and memory storage. We mention
%the recent work by \cite{freitag}. Also the Inexact  Rayleigh Quotient iteration has been 
%recently analyzed by~\cite{xue}.  We finally mention the LOBPCG (Locally Optimal Block Preconditioned Conjugate Gradient Method)
%code developed by  \cite{Kny2001}, which is available under the {\tt hypre} package developed in the
%\cite{LLN01}.

The layout of the rest of the paper is the following. In Section \ref{laplacian} we describe some applications of the lowest eigenpairs of the Laplacian matrix of a graph. The compared state-of-the-art algorithms for eigenpair computation of large and sparse matrices are reviewed in Section \ref{state-of-the-art}. Section \ref{experiments} is devoted to the discussion of the outcomes of the comparison among the algorithms when they run on real network data. We draw our conclusions in Section \ref{conclusion}. 

\section{Why computing some eigenpairs of the graph Laplacian} \label{laplacian}

Let $\mathcal{G}=(V,E,w)$ be a simple (no multiple edges, no self-loops) undirected
weighted graph with $V$ the set of nodes, $|V|=n$, $E$ the set of edges, $|E|=m$, and
$w$ a vector such that $w_k > 0$ is the positive weight of edge $k$, for $k = 1,
\ldots, m$.   The weighted Laplacian of $\mathcal{G}$ is the symmetric matrix
$$G=D-A,$$ where $A$ is the weighted adjacency matrix of the graph and $D$ is the
diagonal matrix of the generalized degrees (the sum of the weights of the incident
arcs) of the nodes. In the following we focus on the spectral properties of the graph Laplacian matrix and on their practical importance.

If $e$ denotes a vector of ones  by definition $De=Ae$ so that $Ge=0$. Thus $e$ is an
eigenvector of $G$ associated to the eigenvalue $\lambda_1=0$. In addition if $x\in
\mathbb{R}^n$ then \begin{equation}x^TGx=\frac12 \sum_{i,j=1}^n
A_{i,j}(x_i-x_j)^2.\label{quadratic}\end{equation} This implies that $G$, besides
symmetric, is positive semidefinite, and hence it has real and nonnegative eigenvalues
that is useful to order $0=\lambda_1\le\lambda_2\le\ldots\le \lambda_n$.

A basic result states that the multiplicity of $0$ as an eigenvalue of $G$ coincides
with the number of the connected components of $\mathcal{G}$. Hence $\lambda_2>0$ if
and only if $\mathcal{G}$ is connected. Fiedler \citet{F73} was one of the pioneers of
the study of the relations between the spectral and the connectivity properties of
$\mathcal{G}$, and for this reason $\lambda_2$  is called Fiedler value or algebraic
connectivity.

Since $G$ is symmetric it admits the spectral decomposition $$G=V\Lambda V^T$$ where
$\Lambda$ is the diagonal matrix such that $\Lambda(i,i) = \lambda_i$ and $V$ is orthogonal, i.e.
$VV^T=I=V^TV$, and its columns are the eigenvectors of $G$. Notice that the
normalization of the eigenvector of $G$ associated with the eigenvalue $\lambda_1=0$
yields $V(:,1)=e/\sqrt{n}$. In addition $V(:,2)$, being the eigenvector associated
with the Fiedler value, will be called Fiedler vector.

Certain applicative problems require the minimization of (\ref{quadratic}) under the
condition that the entries of the vector $x$ belong to some discrete set.

An important example is  discussed in \citet{N10} and concerns graph partitioning. Let
us partition $V$ in two subsets $V_1$ and $V_2$. If we set $x_i=\frac12 (-1)^j$ if
$x_i\in V_j$ then (\ref{quadratic}) is the sum of the weighs of the arcs from one of
the subsets to the other, and is called the cut size. The graph partitioning problem
requires to find $V_1$ and $V_2$ of prescribed dimensions $n_1$ and $n_2$,  in such a
way the cut size is minimized. Actually, all the known methods for finding the minimum
are very demanding, since they reduce to an enumeration of all the ${n \choose
n_1}=\dfrac{n!}{n_1!n_2!}$  possible solutions. However, it is possible to approximate
the minimum by relaxing the constraints on the entries of $x$, allowing them to assume
real values in such a way that $x^Tx=n/4$ and $e^Tx=(n_1-n_2)/2$. Observe that
$$\min_{s.t.\scriptsize{\begin{array}{c}x^Tx=n/4\\
x^Te=(n_1-n_2)/2\end{array}}}\frac12 \sum_{i,j=1}^n
A_{i,j}(x_i-x_j)^2=\min_{s.t.\scriptsize{\begin{array}{c}x^Tx=n/4\\
x^Te=(n_1-n_2)/2\end{array}}}x^TGx.$$ If we set $y=V^Tx$ we obtain
$e^Tx=e^TVy=y_1\sqrt{n}$ $$\min_{s.t.\scriptsize{\begin{array}{c}x^Tx=n/4\\
x^Te=(n_1-n_2)/2\end{array}}}x^TGx= \min_{s.t.\scriptsize{\begin{array}{c}y^Ty=n/4\\
y_1=(n_1-n_2)/(2\sqrt{n})\end{array}}}y^T\Lambda y.$$ Thus, if presented in this
spectral form the problem greatly simplifies. It is easy to find that the minimum is
$\frac{n_1n_2}{n}\lambda_2$ and is obtained when $y_1=(n_1-n_2)/(2\sqrt{n})$,
$y_2=\sqrt{(n_1n_2)/n}$ and $y_i=0$ for $i>2$. Hence, the minimum of the original
problem is obtained for $$x=\frac{n_1-n_2}{2n}e+\sqrt{\frac{n_1n_2}{n}}V(:,2),$$
showing the central role played by the Fiedler vector in the problem.

A second example of the same nature is discussed in \citet{HS03}. In the case where
all the weights are equal to one, the minimization of (\ref{quadratic}), with the
constraint that the entries of $x$ belong to the $n!$ possible permutations of the
integers from $1$ to $n$, allows to find an ordering of the nodes of $\mathcal{G}$ 
that concentrates the entries of $A$ near the main diagonal. For this reason is known
as profile reducing ordering. By relaxing the problem we find as an approximate
solution the ordering induced by the entries of the Fiedler vector. This kind of
applications do not require an accurate computation of the entries of the vector.

A different but equally important application concerns the problem of the computation
of betweenness centrality \citep{N10}. This centrality index quantifies the quantity
of information that passes through a node in order to transit between others.
Actually, for the computation of betweenness centralities, a linear system in $G$ for
every couple of nodes of the network has to be solved. This is actually  equivalent to
the computation of $G^+$, the Moore-Penrose generalized inverse of $G$ \citep{GBS08,
BIG03}. It turns out that $$G^+=\sum_{i=2}^n\frac{1}{\lambda_i}V(:,i)V(:,i)^T.$$ The
use of approximations of $G^+$ obtained by partial sums
$$T^{(k)}=\sum_{i=2}^k\frac{1}{\lambda_i}V(:,i)V(:,i)^T,$$ has been proposed in
\citet{BF2012}. Clearly this implies the computation a certain number of the smallest
eigenpairs of the Laplacian. Moreover, if the eigenvalues $\lambda_{i}$, for
$i=k+1,\ldots,n$ are close to each other it is possible to approximate them by means of a
suitable constant $\sigma$ (for example $\sigma=(\lambda_k+\lambda_n)/2$, or simply
$\sigma=\lambda_k$). In \citet{BF2012} it has been shown that the use of
\begin{eqnarray*} S^{(k)}&=&T^{(k)}+\sum_{j=k+1}^n\frac{1}{\sigma}V(:,j)V(:,j)^T\\
&=&\frac{1}{\sigma}I-\frac{1}{\sigma}V(:,1)V(:,1)^T+\sum_{j=2}^k\bigl(
\frac{1}{\lambda_j}-\frac{1}{\sigma}\bigr)V(:,j)V(:,j)^T. \end{eqnarray*} in the place
of $T^{(k)}$ leads to improved approximations of the centralities. It is important to
note that in order to use $S^{(k)}$ no additional eigenpairs  with respect to
$T^{(k)}$ are requested.

\section{State-of-the-art algorithms for eigenvalue computations of large matrices} \label{state-of-the-art}

Starting from the subspace iteration, which is a generalization of the well-known power
method, a number of iterative procedures have been recently developed  to compute a few eigenpairs
of a large and sparse matrix $A$. In the following, we describe three important methods:
%In the context of graphs, we mention the multilevel algorithm implemented in the HSL\_MC73 routine of the HSL
%mathematical software library, that however, only computes the second smallest eigenpair
%\citep{HS03}.

\begin{itemize}
	\item The Implicitly Restarted Lanczos Method (IRLM).
	\item The Jacobi-Davidson method (JD).
	\item The Deflation Accelerated Conjugate Gradient method (DACG).
\end{itemize}

All the methods are iterative, in the sense that they
compute one or more eigenpairs by constructing a sequence of vectors which approximate the exact
solution. They are all based on the following tasks which must be efficiently carried on:
\begin{enumerate}
	\item Computation  of the the product of the matrix $A$ by a vector. This matrix vector product (MVP) has a cost proportional to the number of nonzero entries of $A$.
	\item Computation of a matrix  $M$, known as preconditioner, that approximates $A^{-1}$ in such a way that
	the eigenvalues of 	$I - M A$ are well clustered around $1$ and in addition
		the computations of $Mv$ and $Av$, being $v$ a generic vector, require comparable CPU
		time.
\end{enumerate}

It is important to stress that IRLM and JD are characterized by an inner-outer iteration, where at every outer iteration
a linear system has to be solved. However, while IRLM requires
to solve these linear systems to a high accuracy, which is strictly related to the accuracy requested
for the eigenpairs, for JD inexact solution of inner linear system is sufficient  to achieve overall convergence. On the other hand, DACG does not require any linear system solution.

%On the other hand, one of the most widely used algorithms is the Lanczos method with implicit
%restart, implemented by ARPACK \citet{arpack97}.
%Alternative approaches are Jacobi-Davidson and Deflation Accelerated Conjugate Gradient
%\citet{BergamaschiPutti02}, that seem to be highly competitive with the Lanczos method. In particular in the
%Jacobi-Davidson method it is still needed to solve inner linear systems, but the factorization is
%avoided and substituted by the use of preconditioned iterative Krylov spaces based methods.
%Deflation Accelerated Conjugate Gradient sequentially computes the eigenpairs by minimizing the
%Rayleigh quotient $q(z)= z^T A z / z^Tz$ over the subspace orthogonal to the eigenvectors previously
%computed.

\subsection{Description of Implicitly Restarted Lanczos Method}
The best known method, the Implicitly Restarted Arnoldi Method (IRAM), is implemented within the ARPACK package \citep{arpack97}
and is also available in the most popular scientific computing packages (Matlab \slash R) .
For symmetric positive definite matrices, IRAM  simplifies to IRLM, Implicitly Restarted Lanczos Method
which reduces the computational cost,  by taking advantage of the symmetry of the problem.
%This variant is takes advance of the Implicitly Shifted QR technique
%that is suitable for large scale problems.

The idea of the Lanczos method is to project the coefficient matrix $A$ onto a subspace
generated by an arbitrary initial vector $\fv_1$ and the matrix $A$ itself,
known as Krylov subspace.
In particular, a Krylov subspace  of dimension $m$ is generated by the following set
of independent vectors:
\[ \fv_1, A \fv_1, \ldots A^{m-1} \fv_1.\]
Actually it is convenient to work with an orthogonal counterpart of this basis and to
organize its vectors as columns of a matrix  $V_{m}$. Then, a symmetric and
tridiagonal matrix
\[T_m = \begin{pmatrix}\alpha_1 & \beta_2 &&& \\
\beta_2 & \alpha_2 & \beta_3 & & \\
 & \ddots & \ddots & \ddots &  \\
&& \beta_{m-1} & \alpha_{m-1} & \beta_{m}  \\
&& &  \beta_{m} & \alpha_{m}
\end{pmatrix} \]
can be computed as
\[ T_m = V_m ^T A V_m. \]
It is well known that the largest eigenvalues  of $T_m$, $\lambda_n^{(m)}, \lambda_{n-1}^{(m)}, \cdots$
converge, as the size of the Krylov
subspace $m$ increases, to the largest eigenvalues of $A$: $\lambda_m, \lambda_{m-1}, \cdots$, while the corresponding eigenvectors
of $A$ can be computed from the homologous eigenvectors of $T_m$ by $\fu_i = V_m \fu_i^{(m)}$.
%\begin{figure}
\begin{algorithm}[t]
	\caption{Implicitly Restarted Lanczos Method}
\begin{minipage}{7cm}
\label{Lanczos}
{\bf Computation of $T_m$ for $A$}.\\
$ v_1 := \,$  unitary norm initial vector. \\
$  v_0 := 0 \,$  \\
  $ \beta_1 := 0 \,$  \\
    {\bf for} $j = 1,2,\cdots,m\,$  \\
\phantom{ {\bf for} }$   w_j := A v_j \, $  \\
\phantom{ {\bf for} }$  \alpha_j :=  w_j^T  v_j  \, $ \\
\phantom{ {\bf for} }$   w_j := w_j - \alpha_j v_j   - \beta_j v_{j-1} \, $ \\
\phantom{ {\bf for} }$  \beta_{j+1} := \left\| w_j \right\|  \, $ \\
\phantom{ {\bf for} }$   v_{j+1} := w_j / \beta_{j+1}  \, $ \\
    {\bf end for} \\
\end{minipage}
%\end{algorithm}
\begin{minipage}{8cm}
%\begin{algorithm}
\label{lan_smallest}
%	\caption{Lanczos Algorithm for the smallest eigenpairs}
{\bf Computation of $T_m$ for $A^{-1}$}.\\
$ v_1 := \,$  unitary norm initial vector. \\
$  v_0 := 0 \,$  \\
  $ \beta_1 := 0 \,$  \\
    {\bf for} $j = 1,2,\cdots,m\,$  \\
    \phantom{ {\bf for} }$   \text{Solve}\ A w_j  =  v_j \, $  \ \text{using the PCG method}\\
\phantom{ {\bf for} }$  \alpha_j :=  w_j^T v_j  \, $ \\
\phantom{ {\bf for} }$   w_j := w_j - \alpha_j v_j   - \beta_j v_{j-1} \, $ \\
\phantom{ {\bf for} }$  \beta_{j+1} := \left\| w_j \right\|  \, $ \\
\phantom{ {\bf for} }$   v_{j+1} := w_j / \beta_{j+1}  \, $ \\
    {\bf end for} \\
\end{minipage}
\end{algorithm}

Such a convergence in many cases is very fast:
roughly $2 N_{eig} \div 3 N_{eig}$
matrix-vector products are usually enough
to compute a small number $N_{eig}$ of the rightmost eigenpairs
to a satisfactory accuracy. This eigenvalue solver exits whenever
the following test is satisfied:
\[\sum_{k=1}^p \frac{1}{p}\frac{\|A \fu_k - \lambda_k \fu_k \|}{\lambda_k} \le \delta,\]
with $\delta$ a fixed tolerance.

Convergence to the smallest eigenvalues is much slower.  Hence, to compute the leftmost part of the
spectrum,  it is more usual to apply the Lanczos process to the inverse of the coefficient matrix $A^{-1}$.
Since $A$ is expected to be large and sparse, its explicit inversion is not convenient from
both CPU time and storage point of view.
Algorithm \ref{Lanczos}, left code, must then be changed since now $w_j$ is computed as the solution of
the linear system $A w_j = v_j$, as reported in Algorithm \ref{Lanczos}, right code.

%\end{figure}
 As before, a Krylov subspace of size roughly $3 N_{eig}$ is sufficient to have $N_{eig}$ leftmost eigenpairs
 to a high accuracy. The complexity of this algorithm is then $3 N_{eig}$ solutions of linear systems
with $A$ as the coefficient matrix.

\noindent
\subsubsection{Solution of the linear system}
The linear system solution needed at every Lanczos step can be solved either
by a direct method (Cholesky factorization)
or by an iterative method such as the Preconditioned Conjugate Gradient (PCG) method.
The former approach is unviable if the system matrix size is large (say $n > 10^4 \div 10^5$)
due to the excessively dense triangular factor provided by the direct factorization.
In such a case the PCG method should be used with the aid of a preconditioner, which
speeds ups convergence. We choose the best known multi purpose preconditioner: the incomplete
Cholesky factorization with no fill-in. Another advantage of the iterative solution is that
the iterative procedure to solve the inner linear system is usually stopped when the following
test is satisfied:
\[ \frac{\|v_j - A w_j\|}{\|v_j\|} \le \delta_{PCG}\]
where the tolerance $\delta_{PCG}$ can be chosen proportional to the accuracy
required for the eigenvectors.

\subsection{Description of the Jacobi-Davidson method}
To compute the smallest eigenvalue this method considers the minimization of the Rayleigh Quotient
\[ q(x) = \frac{x^T A x}{x^T x}  \]
which can be accomplished by setting its gradient to 0, namely
\begin{equation}
\label{nonlin}
 A x - q(x) x = 0.
\end{equation}

Equation (\ref{nonlin}) is a nonlinear system of equations which can be solved by means of the classical
Newton's method in which the Jacobian of (\ref{nonlin}) (or the Hessian of the Rayleigh Quotient) 
is replaced by a simplified formula:
$J(u) \approx A  - q(u)I$, which is shown to maintain the convergence properties of the Newton's method.
The $k$th iterate of this Newton's method hence reads
\begin{eqnarray}
(A  - q(u_k)I) s_k &=& - (A u_k - q(u_k) u_k) \label{linsis} \\
         u_{k+1}   & = & u_k + s_k.
\end{eqnarray}
In practice solution of the system (\ref{linsis}) is known to produce stagnation in the
Newton process. \citet{vdv} proposed to use a projected Jacobian namely
\begin{eqnarray}
(I - u_k u_k^T) (A  - q(u_k)I) (I - u_k u_k^T)s_k &=& - (A u_k - q(u_k) u_k) \label{JDsis} \\
         u_{k+1}   & = & u_k + s_k
\end{eqnarray}
ensuring that the search direction $s_k$ be orthogonal to $u_k$ to avoid stagnation.

Even this corrected Newton iteration may be slow, especially if a good starting point
is not available.
In \citet{vdv} it is proposed
to perform a Rayleigh-Ritz step at every Newton iteration.
In detail, the Newton iterates are collected as columns of a matrix $V$;
then a very small matrix $H = V^T A V$ is computed. The leftmost eigenvector
$y$ is easily computed and a new vector $u$ is obtained as $u = Vy$.
The main consequence of this procedure is the acceleration of the Newton's method
toward the desired eigenvector.

To compute $\lambda_2$, $\lambda_3$, $\cdots$, the previous scheme can be used provided that
the Jacobian matrix is projected onto a subspace orthogonal to the previously computed
eigenvectors. In detail, if $\lambda_j$ is to be computed, the Newton step reads:
\begin{eqnarray}
(I - Q Q^T) (A  - q(u_k)I) (I - Q Q^T)s_k &=& - (A u_k - q(u_k) u_k) \label{JDsismorethanone} \\
         u_{k+1}   & = & u_k + s_k
\end{eqnarray}
where $Q = \left[ v_1 \ v_2 \ \ldots \  v_{j-1} \  u_k \right]$.
In order to maintain the dimension of matrix $H$ sufficiently small,
two additional parameters are usually introduced.  If the size of matrix $H$ is larger than
$m_{\max}$ then only the last $m_{\min}$ columns of matrix $V$ are kept.

Even more than in the Lanczos process, the solution of linear system
(\ref{JDsis}) must be found using an iterative method. This system is usually solved to a very low accuracy so that
in practice  few iterations (20 $\div$ 30) are sufficient to provide a good search direction
$s_k$.
Moreover, it has been proved in \cite{notay} that linear system (\ref{JDsis}) can be solved
by the PCG method, despite of the fact that the system matrix is not
symmetric positive definite. The resulting algorithm is very fast in computing the smallest eigenvalue provided
that a good preconditioner is available for the matrix $A$ in order to solve efficiently the system
(\ref{JDsis}).

\begin{algorithm}[t]
	\caption{Jacobi-Davidson method}
	\label{JD}
\begin{list}{}{\topsep=0.6em\itemsep=0.1em}
	\item Choose unitary starting vector $\fv$.  Initialize empty matrices $H, V$ and $W$, $k=0$.
	\item {\sc while}  $ k < k_{\max} $ \
	\begin{enumerate}{\topsep=0.2em\itemsep=0.0em}
		\item $ k := k+1$.
		\item Orthogonalize $\fv$ against $V$ via modified Gram-Schmidt.
		\item Normalize $v$. Compute $\fw = A \fv$.
		\item $ \ds H :=  \begin{pmatrix} H & V ^T \fw \\ v^T W & v^T w \end{pmatrix},
			\qquad V := [V | v]
			\qquad W := [W | w] $.
		\item Compute the smallest eigenpair $(\theta, \fy)$ of $H$ (with $\|\fy|\| = 1$).
		\item Compute the vector $ \fu := V \fy$ and the associated residual vector $\fr := A \fu - \theta \fu$.
		\item  {\sc if} ${\|\fr\|} < \varepsilon\,  \left(\theta \|\fu\| \right)$  \  STOP
		\item Solve the linear system
			\[(I - \fu \fu^T) (A - \theta I) (I - \fu \fu^T) \fv  = -(A \fu - \theta \fu), \]
            using the  PCG method.
		\end{enumerate}
	\item{\sc end while}
	\end{list}
\end{algorithm}

\subsubsection{Comments on the algorithm}
The sketch of the Jacobi-Davidson algorithm is reported in Algorithm \ref{JD}.
Step 5 implements the Rayleigh-Ritz projection. It is a crucial step for the convergence
of the algorithm, but  requires small CPU time since it consists in the eigensolution
of the usually very small matrix $H$.

Step 8 is the most relevant one from the viewpoint of computational cost. A good projected preconditioner
should be devised in order to guarantee fast convergence of the PCG method. We used here
as the preconditioner $M = (I - \fu \fu^T) P (I - \fu \fu^T) $, with
$P$ the same incomplete Cholesky factorization employed by IRLM. \\[.4em]

\noindent
For the details of this method we refer to the paper~\cite{vdv},
as well as to successive works by \cite{statopulos,fokslevdv,notay}
who analyze both theoretically and experimentally a number of variants of this well known method.

\subsection{Description of the Deflation Accelerated Conjugate Gradient method}
Instead of minimizing $q(x)$ by Newton's method the nonlinear Conjugate Gradient
method can be employed. Differently from the two methods just described,
this one does not need any linear system solution.
Like the JD method, Deflation Accelerated Conjugate Gradient (DACG) computes
the eigenvalues sequentially, starting from the smallest one \citep{bgp97nlaa,berpinsar98}.
The leftmost eigenpairs are computed sequentially, by minimizing the Rayleigh Quotient
over a subspace orthogonal to the previously computed eigenvectors.
Although not as popular as IRLM and JD, this method, which applies only to symmetric positive definite matrices,
has been proven very efficient in the solution of eigenproblems arising from
discretization of Partial Differential Equations (PDEs) in \cite{BergamaschiPutti02}
DACG also proved very suited to parallel implementation
as documented in \cite{bmp12jam} where an efficient parallel matrix vector product has been employed.

Convergence of DACG is strictly related to the relative separation between consecutive eigenvalues,
namely
\begin{equation}
\label{xi}
 \xi_j = \frac{\lambda_j}{\lambda_{j+1}-\lambda_j}.
\end{equation}
 When two eigenvalues are relatively very close,
DACG convergence may be very slow. Also DACG takes advantage of preconditioning
which, as in the two previous approaches,  can be chosen to be the Incomplete Cholesky  factorization.

\begin{algorithm}[t]
	\caption{Deflation Accelerated Conjugate Gradient method}
	\label{DACG}
\begin{list}{}{\topsep=1.0em\itemsep=0.0em\parsep=0.5em}
\item Choose tolerance $\varepsilon$, set $U = 0$.
\item {\sc do}  $\ j=1,p$
  \begin{list}{\theenumi.}{\usecounter{enumi}%
               \leftmargin=3.0em\itemsep=0.2em\topsep=-0.5em}
  \item Choose ${\fx}_0$
        such that $U^T {\fx}_0=0$;  set $k=0$, $\beta_0 = 0$;
  \item
Find the minimum of the Rayleigh Quotient $ q(x) = \dfrac{x^T A x}{x^T x} $
for every $\fx$ such that $U^T {\fx}_0=0$ by a nonlinear preconditioned
conjugate gradient procedure.
\item
Stop whenever the following test is satisfied:
\[   \frac{\|A{\fx}_{k}-q_{k}{\fx}_k \|}{q(\fx_k)} \le \varepsilon\]
        \item
        Set $\begin{displaystyle}
        \lambda_{j}=q_{k},   \qquad
        {\fu}_{j}={\fx}_{k}/\sqrt{\eta}, \qquad U = [U, \fu_j].
        \end{displaystyle}$
  \end{list}
\item{\sc end do}
\end{list}
\end{algorithm}

\subsubsection{Comments on the algorithm.}
The DACG procedure is described in Algorithm \ref{DACG}. The  PCG minimization of the Rayleigh Quotient (Step 2) is carried out by performing a number
of iterations. The main computational burden of a single iteration  is represented by:
\begin{enumerate}
\item One matrix-vector product.
\item One application of the preconditioner.
\item Orthogonalization of the search direction against the previously computed eigenpairs (columns
of matrix $U$). The cost of this step  is increasing with the number of eigenpairs begin sought.
\end{enumerate}
As common in the iterative methods, the number of iterations can not be known in advance.
However, it is known to be proportional to the reciprocal of  the relative separation $\xi$ between
consecutive eigenvalues (Equation (\ref{xi})).

\section{Numerical results and comparisons} \label{experiments}

In this section we experimentally compare the three previously described solvers in the computation
of some of the leftmost eigenpairs of a number of Laplacian matrices  of graphs arising from the following realistic applications covering all four main categories of real networks, namely biological networks, technological networks, information networks, and social networks:

\begin{enumerate}
\item Matrix \prot represents the Laplacian of the protein-protein interaction network of yeast \citep{HMBO01}. In a protein-protein interaction network the vertices are proteins and two vertices are connected by an undirected edge if the corresponding protein interact. 

\item Matrix \ww a symmetrized snapshot of the structure of the Internet at the level of autonomous systems, reconstructed from BGP tables posted by the University of Oregon Route Views Project. This snapshot was created by Mark Newman and is not previously published.

\item Matrix  \www is the Laplacian of the Web network within nd.edu domain \citep{AJB99}. This network is directed but arc direction has been ignored in order to obtain a symmetric Laplacian.

\item Matrix {\sc dblp} is the Laplacian of a graph describing collaboration network of computer scientists. Nodes are authors and edges are collaborations in published papers, the edge weight is the number of publications shared by the authors \citep{F11}.
\end{enumerate}

In Table \ref{char} we report the number of matrix rows (n), the number of matrix nonzero entries (nnz), the average nonzeros per row (anzr), which account for the sparsity of the matrix, and the ratio $\lambda_{51} / \lambda_2$ (gap), which  indicates how the smallest eigenvalues are separated.
Note that the number of nonzeros is computed as $nnz = n + 2m$ where $m$ is the number of arcs in the graph.
\begin{table}[h!]
	\caption{Main characteristics of the sample matrices: size (n), number of nonzero entries (nnz), average nonzeros per row (anzr) and
	ratio between the largest and smallest computed eigenvalues (gap).}
	\begin{center}
		\begin{tabular}{lrrcr}
	matrix  & n & nnz & anzr & gap \\
	\hline
	\prot           & 1\,453 & 5\,344 & 3.7  & 7.28\\
	    \inte          & 22\,963 & 119\,835 & 5.2 & 4.39 \\
	    \www           & 325\,729 & 2\,505\,945 & 7.8 & 23.25 \\
    \dblp           & 928\,498 & 8\,628\,378 & 9.3 & 2.11\\
	\hline
\end{tabular}
\end{center}
\label{char}
\end{table}

%\begin{comment}
\noindent
\begin{figure}[t]
	\centerline{\includegraphics[width=10cm]{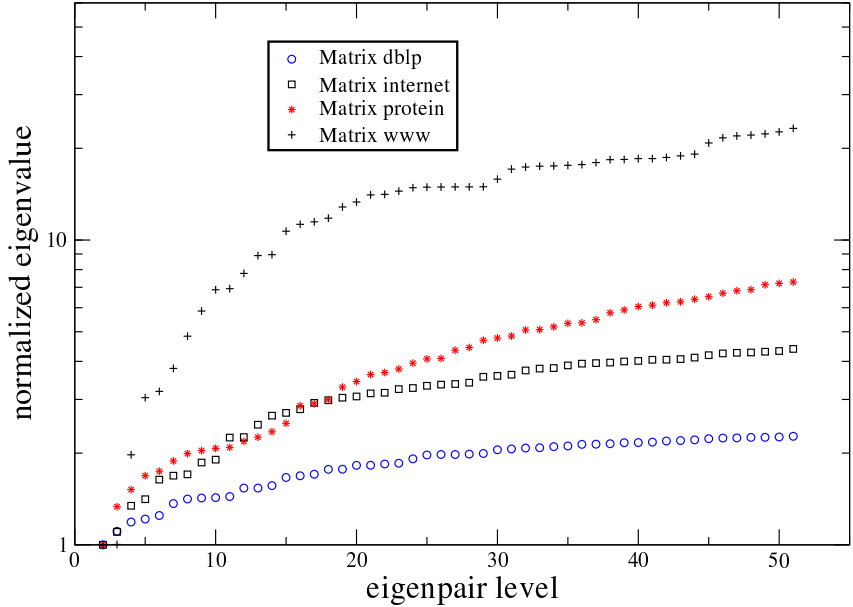}}
\caption{Semilog plot of the distribution of the 50 smallest normalized eigenvalues $(\lambda_j / \lambda_2)$ of the four test matrices. }
\label{distribution}
\end{figure}
%\end{comment}

We also report the distribution of the first 50
normalized eigenvalues for the four test problems in Figure \ref{distribution}.
As mentioned before, a pronounced  relative separation  between consecutive eigenpairs may suggest fast convergence
of the iterative procedures, and particularly so for the DACG method. We notice from Figure \ref{distribution} that for three problems
out of four, with the exception of matrix  \www, the eigenvalues are clustered, thus suggesting a slow convergence
of the iterative solvers. This fact is also accounted for by the ratios $\lambda_{51}/\lambda_{2}$ provided by Table \ref{char}.
The smallest this ratio, the slowest the convergence to the desired eigenvalues.

%\begin{comment}
In the JD implementation two parameters are crucial for its efficiency namely $m_{\min}$ and $m_{\max}$, the smallest and the largest dimension of the subspace where the Rayleigh Ritz projection takes place. After some attempts, we found that $m_{\min} = 5$ and $m_{\max} = 10$ were on the average the optimal values  of such parameters. As for the solution of the Newton linear systems we choose  $ITMAX = 20$ and  we use as the accuracy for the inner linear solver $\delta_{PCG} = 10^{-2}$ .
Regarding IRLM parameters, we set $\delta_{PCG} = 10^{-2} \times \delta$, since the iterative solution
of the inner linear system must be  run to a higher accuracy than that required for the eigenpairs.
The dimension of the Krylov subspace {\tt ncv} for the restarted Lanczos iteration has been chosen as
{\tt ncv} $= 15, 30, 60, 120$, for $N_{eig} = 1, 5, 20, 50$, respectively.
%\end{comment}
The previously described parameters regard memory storage and efficiency for both JD and IRLM. JD is usually 
less demanding than IRLM in terms of memory storage.
If $N_{eig}$ eigenvectors are to be computed, The Jacobi-Davidson method
requires saving of at most $N_{eig} + m_{\max} = N_{eig} + 10$ dense
vectors while {\tt ncv} dense vectors are needed by IRLM. Also the fact that the inner linear system has to be solved
much more accurately by IRLM is accounted for by the choice of parameter $\delta_{PCG}$.

The three solvers,  IRLM, DACG and JD have been  preconditioned by $K^{-1} = \left(L L\t\right)^{-1}$
being $L$ the lower triangular factor of the Cholesky factorization  of $A$ with no fill-in. We reported the results in computing the smallest strictly positive  $N_{eig} = 1, 5, 20, 50$  eigenvalues
with tolerances $\delta = 10^{-3}, 10^{-6}$ for the relative residual using as the exit test:
\[ \frac{\|A \fu_k -  \theta_k \fu_k\|}{\theta_k} < \delta\]
where $\theta_k = \fu_k\t A \fu_k $ is the approximation of the eigenvalue begin $\fu_k$ normalized.

The results regarding the four test problems are summarized in Tables
\ref{protein}, \ref{internet}, \ref{www} and \ref{dblp}, respectively, where  we report CPU times and number of matrix vector products (MVP) for the three codes. The number of linear system solutions of IRLM and JD are also provided (outer iterations). Notice that DACG does not need to solve any linear system. The Fortran implementations of the three solvers have been run on an IBM Power6 at 4.7 GHz and with up to 64 Gb of RAM. The CPU times are expressed in seconds.

\begin{table}[h!]
\caption{Number of linear system solutions (outer its), number of matrix vector products (MVP), and CPU times for DACG, JD, and IRLM on matrix \prot for the computation of the 50 smallest eigenvalues with
two different accuracies ($\delta$).}
\label{protein}
\begin{center}
\begin{tabular}{|c|rr|rrr|rrr|}
        \hline
	&  \multicolumn{2}{|c|}{DACG} & \multicolumn{3}{c|}{JD}   & \multicolumn{3}{|c|}{IRLM} \\
 $\delta$ &           MVP  & CPU & outer its & MVP & CPU  & outer its & MVP & CPU  \\
	\hline
 $10^{-3}$  &3623  &0.46 & 191 & 2501 & {\bf 0.44} & 155 & 2403 & 0.45 \\
 $10^{-6}$  &6513  &0.82 & 293 & 4084 & {\bf 0.68} & 175 & 4381 & 0.72 \\

	\hline
\end{tabular}
\end{center}
\end{table}
Regarding the small-size \prot we only report the results of the computation of $50$ eigenpairs since
computing $N_{eig}=1, 5$ and $20$  eigenvalues is done by every solver in  an almost negligible CPU time  on our
computer.
Table \ref{protein} shows a very similar behavior of the three solvers both in terms of
MVPs and CPU time.

\begin{table}[h!]
\caption{Number of linear system solutions (outer its), number of matrix vector products (MVP), and CPU times for DACG, JD, and IRLM on matrix \inte for the computation of the smallest $1, 5, 20$ and $50$  eigenvalues with two different accuracies ($\delta$).}
\label{internet}
\begin{center}
\begin{tabular}{|cc|rr|rrr|rrr|}
        \hline
	&&  \multicolumn{2}{|c|}{DACG} & \multicolumn{3}{c|}{JD}   & \multicolumn{3}{|c|}{IRLM} \\
        \hline
	$N_{eig}$ & $\delta$ &           MVP  & CPU & outer its & MVP & CPU  & outer its & MVP & CPU  \\
	\hline
	1 & $10^{-3}$  &54 & {\bf 0.48} & 13 & 162  & 0.66 & 21 & 939 & 2.23\\
	5 & $10^{-3}$  &324 & {\bf 0.97} & 33 & 337 & 1.04 & 29 & 1176 & 2.86 \\
20 & $10^{-3}$  &1394  & 3.12 & 82 & 1030 & {\bf 2.89} & 82 & 2231 & 7.41 \\
50 & $10^{-3}$  &4090  &10.66 & 197 & 2383 & {\bf 7.71} & 202 & 5258 &19.29 \\
        \hline
	1 & $10^{-6}$  &93 & {\bf 0.50} & 15 & 159 & 0.67 & 27 & 1927 & 4.40 \\
	5 & $10^{-6}$  &589 & 1.42 & 32 & 454 & {\bf 1.26} & 40 & 2904 & 6.48 \\
20 & $10^{-6}$  &2701  & 5.41 & 116 & 1585 & {\bf 3.98} & 96 & 6724 &15.72 \\
50 & $10^{-6}$  &7591  &18.97 & 285 & 4266 & {\bf 12.93} &211 &14578 &36.31 \\
	\hline
\end{tabular}
\end{center}
\end{table}

\begin{table}[h!]
\caption{Number of linear system solutions (outer its), number of matrix vector products (MVP), and CPU times for DACG, JD, and IRLM on matrix \www for the computation of the smallest $1, 5, 20$ and $50$  eigenvalues with two different accuracies ($\delta$).}
\label{www}
\begin{center}
\begin{tabular}{|cc|rr|rrr|rrr|}
        \hline
	&&  \multicolumn{2}{|c|}{DACG} & \multicolumn{3}{c|}{JD}   & \multicolumn{3}{|c|}{IRLM} \\
        \hline
	$N_{eig}$ & $\delta$ &           MVP  & CPU & outer its & MVP & CPU  & outer its & MVP & CPU  \\
	\hline
	1 & $10^{-3}$  &1262 &     {\bf 50}     &  34 & 1631 & 90   & 37 &18999 & 655 \\
	5 & $10^{-3}$  &5835 & 246 & 63 & 2958 & {\bf 169} & 49 & 24023& 870  \\
20 & $10^{-3}$  &23733  & 1211 &103 &9851 &{\bf 604}  & 81 & 36483 & 1231 \\
50 & $10^{-3}$  &64197  & 4086 & 504 & 22993 &   {\bf 1757}    & 120   & 53742    & 2289    \\
        \hline
	1 & $10^{-6}$  &2121&{\bf 96}&44 &2143&143& 40 &2502 &860 \\
	5 & $10^{-6}$  &9058&428 &109 &5359 &{\bf 358} & 52 & 31244 & 1147\\
20 & $10^{-6}$  &38544& 2329 &373 &18317 &{\bf 1411} & 97 &56578 & 2217 \\
50 & $10^{-6}$  &143102&10935& 986 &44469 &{\bf 3545} &150 &89675 &3853 \\
	\hline
\end{tabular}
\end{center}
\end{table}

\begin{table}[h!]
\caption{Number of linear system solutions (outer its), number of matrix vector products (MVP), and CPU times for DACG, JD, and IRLM on matrix \dblp for the computation of the smallest $1, 5, 20$ and $50$  eigenvalues with two different accuracies ($\delta$).}
\label{dblp}
\begin{center}
\begin{tabular}{|cc|rr|rrr|rrr|}
        \hline
	&&  \multicolumn{2}{|c|}{DACG} & \multicolumn{3}{c|}{JD}   & \multicolumn{3}{|c|}{IRLM} \\
        \hline
	$N_{eig}$ & $\delta$ &           MVP  & CPU & outer its & MVP & CPU  & outer its & MVP & CPU  \\
        \hline
	1 & $10^{-3}$  &178 & {\bf 79}       & 10 & 148 &  82  & 17 & 1359& 421\\
 5 & $10^{-3}$  &797 &  286 & 31 & 450 & {\bf 206} & 41 & 2979 & 937 \\
20 & $10^{-3}$  &3808  &1324  & 107 & 1675 & {\bf 745}   & 103 & 8001  &2447 \\
50 & $10^{-3}$  &11239  &5231  & 315 & 5384 & {\bf 2463}   & 225 & 11111  &4033 \\
        \hline
	1 & $10^{-6}$  &1000 & 319       & 16 & 244 & {\bf 118}  & 24 & 3373 & 1021\\
	5 & $10^{-6}$  &2142 &  690 & 48 & 803 & {\bf 340} & 43 & 5975 & 1813 \\
20 & $10^{-6}$  &7810  &2709  & 171 & 2981 & {\bf 1257}   & 110 & 15628 &4728 \\
50 & $10^{-6}$  &24978  &10571  & 543 & 9654 & {\bf 4665}   & 270 & 34445 &11396 \\
	\hline
\end{tabular}
\end{center}
\end{table}

Analyzing the results in Tables \ref{internet}, \ref{www}, and \ref{dblp} we can make the following observations:
\begin{enumerate}
	\item The IRLM is almost always slower than the remaining two.
	This occurs since  it requires many accurate inner linear system solutions. Actually IRLM
	implicitly computes the largest eigenpairs of $A^{-1}$. The latter matrix is not explicitly formed
	since it would have an excessive number of nonzeros. Only the action of $A^{-1}$ upon a vector
	is computed as a linear system solution. Solving this system to a low accuracy would mean compute
	eigenpairs of a matrix different from $A$ thus introducing unacceptable errors.

\item The JD algorithm displays the best performance in terms of number of MVP and CPU time.
In particular on the largest problem, it neatly outperforms  both DACG and IRLM. The JD method inherits the nice convergence
properties of the Newton's method enhanced by the Rayleigh-Ritz acceleration. Moreover,
it allows very inaccurate (and hence very cheap) solution of the inner linear system.
\item The DACG algorithm provides comparable performances with JD for a very small number of eigenpairs (up to 5)
	and particularly when the eigenvalues are needed to a low accuracy. When the number of eigenvalues
	is large, the reorthogonalization cost prevails and makes this algorithm not competitive.
	For the sample tests presented in this paper, the DACG method is also penalized by the clustering of eigenvalues
	which results in a very small relative separation between consecutive eigenvalues.
	This argument applies also to matrix \www where, apart of the first 10 eigenvalues which are relatively
	well separated, see Figure \ref{distribution}, the remaining ones are as clustered as those of the other test problems.
\item
	When few eigenpairs are to be computed (and hence the reorthogonalization cost is not prevailing)
	JD does not seem particularly sensitive to eigenvalue accuracy. This is not surprising
	as it is based on a Newton iteration. This process is known to converge very rapidly
	in a neighborhood of the solution. For this reason, the transition between
	$\delta = 10^{-3}$ and $10^{-6}$ tolerance is very fast.
\item Despite of the favorable distribution of the leftmost part of its eigenspectrum (see Figure \ref{distribution}), 
	the number of iterations to eigensolve
	matrix  \www is high for all the three solvers. For this test problem, the incomplete
Cholesky factorization with no fill-in preconditioner
	does not provide a satisfactory acceleration.
 The choice of a suitable preconditioner is crucial for the convergence of all iterative methods.
	Devising a more ``dense'' preconditioner 
	would improve the performance of all the   methods
	described and particularly so for IRLM and JD that explicitly require a linear system solution.
\end{enumerate}

%\subsection{Computation of eigenpairs of $D^{-1/2}AD^{-1/2}$}
%{\bf Enrico introduci la motivazione per calcolare gli autovalori di $ D^{-1/2}AD^{-1/2}A$.}
%
%\begin{minipage}{7cm}
%%\begin{figure}[h!]
%	        \centerline{\includegraphics[width=7cm]{w2}}
%		%\caption{Semilog plot of the distribution of the 200 smallest normalized eigenvalues $\left(\dfrac{\lambda_j}{\lambda_1}\right)$
%		%of matrices \inter with and without diagonal scaling. }
%		\label{distribution}
%%	\end{figure}
%\end{minipage}
%\hspace{1cm}
%\begin{minipage}{7cm}
%The eigenvalue distribution of $\tilde{A} = D^{-1/2}AD^{-1/2}A$ and $A$ are obviously different. One may expect
%that the diagonal scaling makes eigenvalues of $\tilde{A}$ more clustered than those of $A$. However,
%we experimentally found that there are only small differences between  the relative eigenvalue
%distribution of the two matrices. As an example we report the distribution of the 200 smallest
%eigenvalues for test case \inte for both the non scaled and scaled matrix.
%\end{minipage}
%\hspace{1cm}
%
%Moreover, the performances of the above iterative solvers,
%which we do not report in this work, are quite similar.

\subsection{Related work}
We selected to use for the three methods the established implementations without making  optimization
to  any of them. However, it is worth mentioning that much work is being devoted particularly  to the 
Arnoldi method (the non symmetric counterpart of the Lanczos Method) in order to reduce its computational 
cost and memory storage. We refer e.g. to the recent work by \cite{freitag}. 

Other methods are efficiently employed for computing a number of eigenpairs of sparse matrices.
Among these, we mention the Rayleigh Quotient iteration whose inexact variant  has been
recently analyzed by~\cite{xue}.  A method which has some common features with DACG is  LOBPCG 
(Locally Optimal Block Preconditioned Conjugate Gradient Method)
which has been  proposed by  \cite{Kny2001}, and is currently available under the {\tt hypre} package developed in the
\cite{LLN01}.

\section{Conclusion} \label{conclusion}

We experimentally compare three iterative state-of-the-art algorithms for computation of eigenpairs of large and sparse matrices: the Implicitly Restarted Lanczos Method, the Jacobi-Davidson method, and the Deflation Accelerated Conjugate Gradient method. We uniformly implemented the algorithms and ran them in order to compute some of the smallest eigenpairs of the Laplacian matrix of real-world networks of different sizes.

The iterative approach followed in this work seems to be particularly suited for the Laplacian matrices
presented since it fully exploits the sparsity of the matrices involved. Each of our realistic test
cases, indeed, has a very high degree of sparsity as accounted for by the very small number of nonzeros per row.

Contrary to what observed for matrices arising from discretization of Partial Differential Equations \citep{BergamaschiPutti02}, 
where especially the smallest eigenvalues are well separated, here the high clustering of lowest eigenvalues is disadvantageous for the 
Deflation Accelerated Conjugate Gradient algorithm. As for the Implicitly Restarted Lanczos Method, the need to solve the inner linear systems 
to a high accuracy makes this method less attractive for large eigenproblems.The Jacobi-Davidson procedure is less sensitive to the 
clustering thanks to the Rayleigh-Ritz projection; moreover, for this method, inexact and hence efficient solution of the inner linear 
systems is sufficient to achieve overall convergence. All in all, the Jacobi-Davidson algorithm is performing the best. 
This might be valuable information for popular scientific computing environments (Matlab \slash R), which only implement 
the Implicitly Restarted Lanczos Method. 

All the proposed algorithms are well-suited to parallelization on supercomputers. The most important kernel is represented
by the matrix-vector product which can be efficiently implemented in parallel environments. 
Also  application of preconditioner, which is one of the most time-consuming task, 
in its turn can be devised as a product of sparse matrices
as e.g. in the ``approximate inverse preconditioner'' approach (see the review article by \cite{benzi}). 
The DACG method has been successfully parallelized as
documented in~\cite{bmp12jam}, however all the iterative solvers described here, being based on the same
linear algebra kernels, could be implemented in parallel with the same satisfactory results

In our implementation we used a general purpose preconditioner, obtained by means of incomplete Cholesky factorization with no fill in. 
Certainly, the three methods would greatly benefit from the use of a more specific preconditioner. This point will be a topic of future research.

\end{document}